\documentclass[11pt]{article}

\usepackage[margin=1in]{geometry}
\usepackage{amsmath,amsthm}
\usepackage{newtxtext}
\usepackage{newtxmath}
\usepackage{bm}
\usepackage{setspace}
\usepackage{mathtools}
\let\mathscr\mathcal
\usepackage{booktabs}
\usepackage{graphicx}
\usepackage{microtype}
\usepackage{enumitem}
\usepackage{float}
\usepackage{placeins}
\usepackage[hidelinks]{hyperref}
\usepackage[nameinlink,noabbrev]{cleveref}
\numberwithin{equation}{section}
\newtheorem{theorem}{Theorem}[section]
\newtheorem{proposition}[theorem]{Proposition}
\newtheorem{lemma}[theorem]{Lemma}
\newtheorem{corollary}[theorem]{Corollary}
\theoremstyle{definition}

\newtheorem{assumption}[theorem]{Assumption}

\newtheorem{remark}[theorem]{Remark}
\crefname{assumption}{Assumption}{Assumptions}

\newcommand{\R}{\mathbb R}
\newcommand{\E}{\mathbb E}

\newcommand{\1}{\mathbf 1}
\newcommand{\norm}[1]{\left\lVert#1\right\rVert}
\begin{document}
\title{Optimal Control of Scalar Diffusions with Finite Action Sets}
\author{Linyu Miao\thanks{School of Mathematical Sciences, Nankai University,
Tianjin 300071, China. Email:
\href{mailto:mly202210@126.com}{mly202210@126.com}.}}
\date{}
\maketitle

\begin{abstract}
We study a scalar finite-horizon controlled diffusion with a finite action
set. An action may affect the drift, the instantaneous variance, or both, and
performance is assessed by a uniformly convex terminal deviation cost. The
model represents systems in which decisions are made continuously in time but
only finitely many operating settings are available. Under uniform
nondegeneracy and an explicit shape condition on the terminal cost, we show
that the value function is classical at positive times, establish global
positive upper and lower bounds for its second state derivative, and prove
that the ratio of its first and second state derivatives is strictly
increasing. This ordering implies that each action's optimality region is
connected, the active actions occur in monotone order, and adjacent regions
are separated by unique switching boundaries. It also yields an exact
coefficient-based screen that removes actions that cannot be optimal on any
open state region. For smooth continuous-action benchmarks, both the value
loss from restricting the action set and the additional loss from rounding
the continuous optimal feedback to the nearest available action decrease
quadratically with the action-set resolution. Numerical experiments illustrate
the resulting action screens, threshold policies, and implementation losses.
\end{abstract}



\section{Introduction}
We study finite-action stochastic control for a scalar diffusion over a
finite horizon. The controller repeatedly selects from finitely many
admissible actions; an action may affect the drift, the instantaneous
variance, or both. This formulation is motivated by systems in which the
decision itself is genuinely discrete. Production and service facilities
offer finitely many operating levels, communication systems select among
prespecified transmission settings, portfolios are rebalanced in integer
lots rather than arbitrary fractions \cite{BertsimasEtAl1999,
BonamiLejeune2009}, and insurance or reinsurance contracts are often offered
in a finite collection of coverage levels. In such settings no separate cost
is needed to make the action set discrete: discreteness is a property of the
available decisions, not a substitute for an otherwise continuous actuator.

The model can look superficially similar to optimal switching because the
selected action changes when the state crosses an action boundary. The
economic primitive is different. Optimal switching describes a persistent
incumbent regime and makes the cost of replacing that regime an essential
part of the decision; its value is indexed by the current regime. Our problem
describes repeated selection of a current action from a finite set. The
previous action is not an additional state variable, and a separate
changeover cost is not required for the model to be meaningful. Thus the two
formulations address complementary operational settings rather than two
notations for the same problem.

This distinction is also visible in the literature. Classical and recent
optimal-switching models analyze technologies, production configurations, or
market positions for which changing the incumbent regime incurs an explicit
transition cost \cite{DjehicheHamadenePopier2009,LvWuZhang2018,
LuSerfozo1984,VandeVate2021,OrmeciVandeVate2024}. Another line of work studies
systems in which the action can be selected directly and no separate
changeover cost is needed. Examples include controlled Brownian motion with
finitely many controls \cite{Rath1977,Sheng1982}, adjustable service rates
\cite{GeorgeHarrison2001,AtaShneorson2006}, dynamic risk taking
\cite{SubramanianYang2017}, and drift--variance choices in new-venture
development \cite{WangZenios2022}. These papers establish threshold or
single-switch structures in important special models. They do not provide a
global answer for a finite, possibly time-varying set of actions that may
jointly alter drift and variance: which actions can ever be optimal, whether
an action can disappear and later reappear as the state varies, and how many
state thresholds are required.

At the level of mathematical analysis, general controlled-diffusion theory
identifies the value function as the unique viscosity solution of the
corresponding Hamilton--Jacobi--Bellman equation under standard comparison
assumptions \cite{Krylov1980,FlemingSoner2006,
CrandallIshiiLions1992}. This well-posedness statement does not determine the
global geometry of the action regions. Nor does the action-quantization
literature, which treats a finite set as an approximation to an underlying
compact action space \cite{SaldiYukselLinder2016,
SaldiLinderYuksel2018,ReisingerForsyth2016}. Here the finite set is the
implemented decision space. A continuous-action model is only a benchmark,
and one must separately ask how much can be lost by rounding its feedback to
the available actions.

Two features make the desired structural result nontrivial. First, the finite
action set is nonconvex. Although dynamic programming readily takes a minimum
over finitely many generators, convexity of the terminal cost does not by
itself guarantee the global derivative ordering needed to rule out
disconnected action regions. Second, when actions affect both drift and
variance, comparisons depend on both the first and second state derivatives
of the value function. The minimizing branch changes across action
boundaries, so the relevant transformed equation has a nonsmooth finite-
envelope flux. Static comparison of the action coefficients alone is
therefore insufficient.

Our first contribution is a global regularity and ordering result. Under
uniform nondegeneracy and an explicit terminal shape condition, we prove
positive upper and lower bounds for the second state derivative and show that
the ratio of the first to the second state derivative is strictly increasing.
The proof uses the inverse of the first derivative, whose evolution is a
uniformly parabolic divergence-form equation, together with a strict spatial
comparison argument that remains valid across the flux junctions. This turns
the static ordering of generator lines into a global ordering of action
regions.

Our second contribution is an exact characterization of the optimal
feedback. At every time, only vertices of the lower convex hull of the
explicit drift--variance action points can be optimal on open state regions.
The retained actions occur in hull order, each retained action has one
connected region, and adjacent regions are separated by a unique state
threshold. The hull calculation is consequently an exact screening rule,
not a heuristic reduction. With time-varying coefficients it applies
pointwise: an action can be screened at some times and become optimal at
others. The drift-only and variance-only models follow as simpler cases.

Our third contribution quantifies the use of a continuous-action planning
model when only finitely many actions can be implemented. We compare the exact
finite-action optimum with the policy obtained by solving the continuous-
action model and rounding its feedback to the nearest available action. For
smooth actuator families, both the loss caused by restricting the available
actions and the additional loss caused by this rounding procedure decay
quadratically with the resolution of the finite action set. Numerical
experiments display the ordered action regions, time-dependent screening, and
the operational gains from solving the finite-action problem directly when
the available actions are sparse.

Section~\ref{sec:model} formulates the control problem and assumptions.
Sections~\ref{sec:inverse}--\ref{sec:hull} establish the derivative bounds,
strict order, and ordered threshold feedback. Section~\ref{sec:scope}
develops exact action screening and the continuous-action comparison, and
Section~\ref{sec:numerics} reports the numerical experiments. Proofs of the
main analytical results are collected in the appendices.

\section{Weak formulation, HJB equation, and structural assumptions}\label{sec:model}
Fix a horizon $T>0$ and a finite action set $\mathcal I=\{1,\dots,N\}$, and
use the time-to-go variable $\tau$.  For an initial pair $(\tau,x)$, let
$\mathfrak A_{\tau,x}$ be the collection of tuples
$(\Omega,\mathcal F,(\mathcal F_s)_{0\le s\le\tau},\mathbb P,W,I,X)$
such that the filtration satisfies the usual conditions, $W$ is a standard
Brownian motion, $I$ is an $\mathcal I$-valued progressively measurable
process, and $X$ is an adapted continuous process satisfying
\begin{equation}\label{eq:state}
 X_s=x-\int_0^s\beta_{I_r}(\tau-r)\,dr
 +\int_0^s\sqrt{a_{I_r}(\tau-r)}\,dW_r,\qquad 0\le s\le\tau.
\end{equation}
The value function is
\begin{equation}\label{eq:value}
 v(\tau,x):=\inf_{\mathfrak A_{\tau,x}}\E\bigl[G(X_\tau)\bigr],
 \qquad v(0,x)=G(x).
\end{equation}
Thus the infimum is over admissible weak systems, not merely over controls
defined on one preassigned probability space.  This formulation is
convenient because the discontinuous Markov selector constructed later need
not generate a unique strong solution.

We use standard parabolic spaces.  If $Q=(t_0,t_1)\times O$, then
$C^{\theta/2,\theta}(Q)$ denotes the parabolic H\"older space and
$H^{-1}(O)$ is the dual of $H^1_0(O)$.  Norms are written
$\norm{\cdot}_{L^p(Q)}$ and $\norm{\cdot}_{\infty,Q}$.  An inequality
$D_{xx}f\ge c$ in distributions means
\begin{equation}\label{eq:distribution-order}
 \int_{\R}f\varphi_{xx}\,dx
 \ge c\int_{\R}\varphi\,dx
\end{equation}
for every nonnegative $\varphi\in C_c^\infty(\R)$.
In one dimension this is equivalent to convexity of
$x\mapsto f(x)-cx^2/2$.

We impose the following assumptions throughout unless stated otherwise.
\begin{assumption}\label{ass:main}
For some $\alpha\in(0,1)$:
\begin{enumerate}
\item $a_i,\beta_i\in C^{\alpha/2}([0,T])$, and there are constants
\[
 0<\underline a\le a_i(\tau)\le\overline a,
 \qquad |\beta_i(\tau)|\le \overline\beta .
\]
\item $G\in C^{3+\alpha}_{\mathrm{loc}}(\R)$, $G''$ is globally H\"older continuous, and
\[
 0<m\le G''(x)\le M<\infty,\qquad x\in\R.
\]
Let $x_0$ be the unique minimizer of $G$.
\item The terminal derivative ratio
\[
 R_G(x):=\frac{G'(x)}{G''(x)}
\]
is strictly increasing.  Quantitatively, we assume
\begin{equation}\label{eq:RG-slope}
 0<r_0\le R_G'(x)\le r_1<\infty,
 \qquad x\in\R.
\end{equation}
\end{enumerate}
\end{assumption}
Uniform ellipticity is structural: it supplies positive-time Bellman regularity and the positive kernels used to turn weak convexity into a strict lower bound for $D_{xx}v$.  The last assumption is stronger than convexity.  Indeed,
\[
 R_G'=1-\frac{G'G'''}{(G'')^2},
\]
so uniform convexity alone does not determine its sign.

The condition is satisfied, for example, whenever $G'G'''\le0$, in which
case $R_G'\ge1$.  One nonquadratic family with this property is
\begin{equation}\label{eq:loss-examples}
 G_1(x)=\frac{x^2}{2}
 +\varepsilon\bigl(\sqrt{1+x^2}-1\bigr),
 \qquad \varepsilon>0.
\end{equation}
Uniform convexity alone is insufficient.  If
$G''(x)=2+\varepsilon\sin x$, $0<\varepsilon<1$, and $G'(0)=0$, then
$G'(x)=2x+\varepsilon(1-\cos x)$, while $R_G'$ is negative at
sufficiently large positive multiples of $2\pi$ and positive at odd
multiples of $\pi$.

\begin{remark}[Common uncontrolled affine state term]\label{rem:affine-reduction}
The model also covers a common uncontrolled affine drift without enlarging
the state.  In calendar time, suppose
\[
 dY_s=\{A(s)Y_s+c(s)-\widehat\beta_{I_s}(s)\}\,ds
 +\sqrt{\widehat a_{I_s}(s)}\,dW_s.
\]
Set
\[
 \phi(s)=\exp\!\left(\int_s^T A(r)\,dr\right),\qquad
 C(s)=\int_s^T\phi(r)c(r)\,dr,
 \qquad Z_s=\phi(s)Y_s+C(s).
\]
Since $\phi'=-A\phi$ and $C'=-\phi c$, It\^o's formula gives
\[
 dZ_s=-\phi(s)\widehat\beta_{I_s}(s)\,ds
 +\phi(s)\sqrt{\widehat a_{I_s}(s)}\,dW_s,
 \qquad Z_T=Y_T.
\]
After reversing time, this is \eqref{eq:state} with
$\beta_i=\phi\widehat\beta_i$ and $a_i=\phi^2\widehat a_i$.
Accordingly, boundedness, H\"older continuity, and uniform nondegeneracy are
preserved under the corresponding assumptions on the original coefficients.
The commonality of $A$ is essential to this reduction; action-dependent state
multipliers are outside the scope of the present model.
\end{remark}

Dynamic programming gives the Bellman equation
\begin{equation}\label{eq:hjb}
 \begin{aligned}
 v_\tau&=\mathcal F(\tau,v_x,v_{xx}),\\
 \mathcal F(\tau,p,r)
 &:=\min_{i\in\mathcal I}
 \left\{-\beta_i(\tau)p+\frac{a_i(\tau)}2r\right\}.
 \end{aligned}
\end{equation}
Standard quadratic-growth comparison identifies \eqref{eq:value} with the
unique viscosity solution of \eqref{eq:hjb}; see
\cite{FlemingSoner2006,CrandallIshiiLions1992}.  At this stage no classical
derivative of $v$ is used.  The structural ellipticity is the elementary
inequality
\begin{equation}\label{eq:uniform-parabolicity}
 \frac{\underline a}{2}\xi
 \le \mathcal F(\tau,p,r+\xi)-\mathcal F(\tau,p,r)
 \le \frac{\overline a}{2}\xi,
 \qquad \xi\ge0.
\end{equation}
The passage from the viscosity solution to a classical positive-time
solution is made only after the distributional estimates in the next
section have been proved.

\section{Derivative bounds, the inverse equation, and strict spatial order}\label{sec:inverse}

\subsection{Distributional second-derivative bounds and positive-time regularity}
The action set is finite, so convexity cannot be proved by averaging two controls.  The following result instead uses the parabolic equation.

\begin{proposition}[Global bounds for the second derivative]\label{prop:second-derivative}
Under \cref{ass:main}, for every $T_0\le T$ there are constants
$0<c_{T_0}\le C_{T_0}<\infty$, depending only on
$\underline a,\overline a,\overline\beta,m,M,\alpha$, and $T_0$, such that
\begin{equation}\label{eq:second-derivative-bounds}
 \begin{gathered}
 c_{T_0}m\le v_{xx}(\tau,x)\le C_{T_0}M,\\
 0<\tau\le T_0,\qquad x\in\R.
 \end{gathered}
\end{equation}
Before positive-time regularity is invoked, \eqref{eq:second-derivative-bounds} is understood as an inequality of distributions.
\end{proposition}

\begin{proof}
See \cref{app:second-derivative}.  The proof is kept there as one continuous
argument: smooth approximation, a maximum-principle proof of convexity,
uniform increment estimates, distributional passage, and finally the
pointwise upgrade.
\end{proof}

\begin{corollary}[Classical positive-time solution]\label{cor:regularity}
For some $\theta\in(0,\alpha)$,
\begin{equation}\label{eq:positive-regularity}
 v\in C^{1+\theta/2,2+\theta}_{\mathrm{loc}}
 ((0,T]\times\R).
\end{equation}
On every open set on which one action is the unique minimizer, $v$ has the
higher regularity of the corresponding linear equation.  The bounds
\eqref{eq:second-derivative-bounds} hold pointwise for $\tau>0$.
\end{corollary}
\begin{proof}
The approximation constructed in \cref{app:second-derivative} satisfies
\eqref{eq:uniform-parabolicity} with the same constants and has locally
uniformly bounded quadratic growth.  The concave parabolic Bellman estimates
of \cite{Krylov1980,Wang1992}, followed by stability, give
\eqref{eq:positive-regularity}.  The distributional density obtained in
\cref{prop:second-derivative} must therefore equal the classical $v_{xx}$.
Inside a unique-action region the equation is linear, so the usual interior
Schauder estimates apply \cite{LadyzhenskayaEtAl1968,Lieberman1996}.
\end{proof}

We shall also need a bound for the inverse of $v_x$ at infinity.
\begin{lemma}[Bounded inverse displacement]\label{lem:displacement}
Let $I_0=(G')^{-1}$.  For every $T_0\le T$ there is $K_{T_0}$ such that
\begin{equation}\label{eq:marginal-displacement}
 \begin{aligned}
 |v_x(\tau,x)-G'(x)|&\le K_{T_0},\\
 |(v_x(\tau,\cdot))^{-1}(p)-I_0(p)|
 &\le K_{T_0}/m .
 \end{aligned}
\end{equation}
\end{lemma}
\begin{proof}
See the final part of \cref{app:second-derivative}, where both inequalities are
derived directly from the control representation.
\end{proof}

\subsection{The inverse of the first derivative}
By \cref{prop:second-derivative,lem:displacement}, $p=v_x(\tau,\cdot)$ is an increasing bijection of $\R$.  Let
\begin{equation}\label{eq:X-def}
 \Xi(\tau,p):=(v_x(\tau,\cdot))^{-1}(p),
 \qquad \Xi(0,p)=I_0(p).
\end{equation}
This already gives existence of the transformed function; no independent
existence construction for \(\Xi\) is required.  For positive time, the
inverse-function theorem and \cref{eq:positive-regularity} give
\begin{equation}\label{eq:Xi-p-basic}
 \Xi_p(\tau,p)=\frac{1}{v_{xx}(\tau,\Xi(\tau,p))}>0,
\end{equation}
with local H\"older continuity.  The purpose of the equation below is to
identify the evolution inherited by \(\Xi\) and to prove uniqueness in a
natural weak class, which is needed for the comparison argument in
\cref{sec:strict}.
If $H(\tau,p)=\sup_x\{px-v(\tau,x)\}$, then $\Xi=H_p$ and
\begin{equation}\label{eq:H-dual}
 \begin{aligned}
 H_\tau&=\mathcal M(\tau,p,H_{pp}),\\
 \mathcal M(\tau,p,r)
 &:=\max_i
 \left\{\beta_i(\tau)p-\frac{a_i(\tau)}{2r}\right\}.
 \end{aligned}
\end{equation}
Differentiating in $p$ gives the global divergence equation
\begin{equation}\label{eq:X-PDE}
 \Xi_\tau=\partial_p\mathcal M(\tau,p,\Xi_p).
\end{equation}

\begin{theorem}[Inherited inverse equation and uniqueness]\label{thm:inverse}
The function \(\Xi\) defined in \eqref{eq:X-def} satisfies
\eqref{eq:X-PDE} and is unique in the whole-line weak class specified in
\cref{app:inverse}.  Moreover,
\begin{equation}\label{eq:X-slope}
 \begin{gathered}
 \frac1{C_{T_0}M}\le \Xi_p\le\frac1{c_{T_0}m},\\
 \norm{\Xi-I_0}_{\infty,(0,T_0)\times\R}
 \le K_{T_0}/m.
 \end{gathered}
\end{equation}
On every compact positive-time cylinder it has the regularity supplied by
uniformly parabolic one-dimensional divergence equations.  Its gradient has
a jointly continuous representative, including at flux junctions and at the
initial slice on compact sets; second derivatives are asserted only inside
open unique-action regions.
\end{theorem}

\begin{proof}
Existence and the slope bounds are inherited from
\cref{eq:X-def,eq:Xi-p-basic,eq:second-derivative-bounds,eq:marginal-displacement}.
The Legendre identity gives \eqref{eq:X-PDE}, first for the smooth Bellman
approximations and then weakly by stability.  The whole-line uniqueness
calculation and the initial-slice continuity argument are given in
\cref{app:inverse}.
\end{proof}

\begin{lemma}[Joint continuity of $\Xi_p$]\label{lem:Xp-continuity}
For every $T_0\le T$ and $P<\infty$, $\Xi_p$ has a positive jointly
continuous representative on $[0,T_0]\times[-P,P]$.  In particular,
\begin{equation}\label{eq:Xp-uniform-zero}
 \lim_{q\downarrow0}\ \sup_{0\le\tau\le T_0}
 |\Xi_p(\tau,\pm q)-\Xi_p(\tau,0)|=0.
\end{equation}
\end{lemma}
\begin{proof}
See \cref{app:inverse}.  The proof checks the Bellman regularity hypotheses
on compact cylinders and then uses
$\Xi_p=1/v_{xx}(\tau,\Xi)$; it does not infer gradient continuity from weak
continuity of $\Xi$.
\end{proof}

\subsection{From the inverse equation to strict spatial order}\label{sec:strict}
Define
\[
 R(\tau,x)=\frac{v_x(\tau,x)}{v_{xx}(\tau,x)}.
\]
\begin{theorem}[Strict spatial order]\label{thm:strict-order}
Under \cref{ass:main}, for every \(\tau>0\) the map
\(x\mapsto R(\tau,x)\) is continuous, strictly increasing, and maps
\(\R\) onto \(\R\).
\end{theorem}
\begin{proof}
See \cref{app:strict-order}.  The auxiliary functions used there are
explicit transforms of the already defined \(\Xi\); their existence and
local H\"older regularity are inherited from
\cref{thm:inverse,lem:Xp-continuity}, rather than obtained from a second
existence construction.
\end{proof}

This result is the dynamic link between the terminal-loss assumption and
action screening.  Every candidate generator slope is crossed once and only
once, so the lower-hull order becomes an ordered list of state regions rather
than a collection of disconnected active components.  The next section turns
this analytical statement into the implementable threshold feedback.

\long\def\StrictOrderAppendix{%
\section{Comparison and strict positivity for the derivative ratio}
\label{app:strict-order}
\subsection{Equations on the two signs of the first derivative}
The logarithmic coordinates are taken separately on the two signs of $p=v_x$.  For $p=e^y>0$, set
\begin{equation}\label{eq:rho-plus}
 \begin{aligned}
 \rho^+(\tau,y)&:=e^y\Xi_p(\tau,e^y),\\
 \Phi^+(\tau,r)&:=\max_i
 \left\{\beta_i(\tau)-\frac{a_i(\tau)}{2r}\right\}.
 \end{aligned}
\end{equation}
For $p=-e^y<0$, set
\begin{equation}\label{eq:rho-minus}
 \begin{aligned}
 \rho^-(\tau,y)&:=e^y\Xi_p(\tau,-e^y),\\
 \Phi^-(\tau,r)&:=\max_i
 \left\{-\beta_i(\tau)-\frac{a_i(\tau)}{2r}\right\}.
 \end{aligned}
\end{equation}
Then
\begin{equation}\label{eq:rho-PDE}
 \rho^\pm_\tau=(\Phi^\pm(\tau,\rho^\pm))_{yy}
 +(\Phi^\pm(\tau,\rho^\pm))_y.
\end{equation}
Each flux is continuous, strictly increasing, and locally Lipschitz on $(0,\infty)$; at an action junction its derivative may have a corner.  From \eqref{eq:X-slope},
\begin{equation}\label{eq:rho-barriers}
 \lambda_{T_0}e^y\le\rho^\pm(\tau,y)
 \le\Lambda_{T_0}e^y.
\end{equation}
The initial functions are
 \[
 \begin{aligned}
 \rho^+(0,y)&=R_G(I_0(e^y)),\\
 \rho^-(0,y)&=-R_G(I_0(-e^y)).
 \end{aligned}
\]
and \eqref{eq:RG-slope} makes both strictly increasing.

For the proof it is preferable to return from $y$ to the magnitude
$q=e^y$.  Define
\begin{equation}\label{eq:F-branches}
 \begin{aligned}
 F^+(\tau,q)&=q\Xi_p(\tau,q),\\
 F^-(\tau,q)&=q\Xi_p(\tau,-q),
 \qquad q>0.
 \end{aligned}
\end{equation}
Then $F^\pm(\tau,q)=\rho^\pm(\tau,\log q)$ and
\begin{equation}\label{eq:F-PDE}
 \begin{aligned}
 F^\pm_\tau
 &=\partial_q\left[q^2\partial_q
 \Phi^\pm(\tau,F^\pm)\right],\\
 \lambda_{T_0}q&\le F^\pm(\tau,q)
 \le\Lambda_{T_0}q.
 \end{aligned}
\end{equation}
The factor $q^2$ cancels the apparent singularity of
$\Phi_r(\tau,F)\asymp F^{-2}$ at $q=0$.

\begin{lemma}[Local Kato inequality]\label{lem:local-kato}
Let $\widetilde Q$ be an open cylinder and let $Q\Subset\widetilde Q$.
Suppose $F_1,F_2$ are bounded weak solutions of \eqref{eq:F-PDE} on
$\widetilde Q$ and satisfy the local slope bounds in \eqref{eq:F-PDE}. Put
\[
 d=F_1-F_2,
 \qquad e=\Phi(\tau,F_1)-\Phi(\tau,F_2).
\]
Then $d$ and $e$ have the same sign and, in distributions on $Q$,
\begin{equation}\label{eq:local-kato}
 \partial_\tau d^-
 \le \partial_q\{q^2\partial_q e^-\}.
\end{equation}
The statement remains valid when $\Phi$ is the finite maximum in
\eqref{eq:rho-plus} or \eqref{eq:rho-minus}.
\end{lemma}
\begin{proof}
Choose $Q'$ with $Q\Subset Q'\Subset\widetilde Q$.  On the compact range of
$F_1,F_2$ there are constants $0<\ell\le L$ such that every branch of
$\Phi$ has slope in $[\ell,L]$.  Smooth the branches without changing these
bounds, denote the resulting approximants by $\phi_j^n$, and set
\begin{equation}\label{eq:smooth-max-flux}
 \Phi_n(\tau,r)=\frac1n\log\!\left(
  \sum_j e^{n\phi_j^n(\tau,r)}\right).
\end{equation}
Then $\Phi_n\to\Phi$ locally uniformly and
$\ell\le\partial_r\Phi_n\le L$.  On $Q'$, the factor $q^2$ is bounded above
and away from zero.  Smooth Dirichlet approximations therefore converge
strongly in $L^2(Q)$, while their flux gradients converge weakly in $L^2(Q)$;
strong monotonicity identifies the limit by Minty's argument.  These are the
standard local stability statements for monotone quasilinear parabolic
equations \cite{AltLuckhaus1983,Carrillo1999}.  Thus it remains only to
record the truncation calculation that determines the sign.

Put $d_n=F_1^n-F_2^n$ and
$e_n=\Phi_n(\tau,F_1^n)-\Phi_n(\tau,F_2^n)$.  Let $S_\delta$ be a smooth
nonincreasing approximation of $\1_{(-\infty,0)}$.  Since $\Phi_n$ is
strictly increasing, $d_n$ and $e_n$ have the same sign.  Multiplication of
$(d_n)_\tau=\partial_q(q^2(e_n)_q)$ by $-S_\delta(e_n)$ gives
\begin{align}
 -S_\delta(e_n)(d_n)_\tau
 &=-\partial_q\{q^2S_\delta(e_n)(e_n)_q\}
   +q^2S_\delta'(e_n)|(e_n)_q|^2
 \label{eq:kato-defect}\\
 &\le-\partial_q\{q^2S_\delta(e_n)(e_n)_q\},
 \nonumber
\end{align}
because $S_\delta'\le0$.  Against a nonnegative test function compactly
supported in $Q$, first let $\delta\downarrow0$.  This gives
\[
 \partial_\tau(d_n)^-
 \le \partial_q\{q^2\partial_q(e_n)^-\}
 \quad\text{in }\mathcal D'(Q).
\]
There are no endpoint terms because the test is compactly supported.  First
let $\delta\downarrow0$ for the smooth problems.  Strong $L^2$ convergence of
the solutions and fluxes, together with the uniform local $H^1$ estimate for
the truncated flux difference, then permits $n\to\infty$ and proves
\eqref{eq:local-kato}.  This is the local version of the standard Kato
contraction argument \cite{Otto1996,Carrillo1999}; the displayed calculation
also shows why the finite-envelope kink causes no extra defect measure.
\end{proof}

\begin{lemma}[Interior positivity for the forward adjoint equation]
\label{lem:adjoint-positivity}
Let $D=(q_-,q_+)\Subset(0,\infty)$ and let
$d,b$ be Borel functions on $(0,t_*)\times D$ satisfying
\[
 0<d_0\le d\le d_1,
 \qquad |b|\le b_1.
\]
Suppose
\[
 w\in C((0,t_*]\times D)\cap C([0,t_*];L^1(D)),
 \qquad w\ge0,
\]
is a distributional supersolution of
\begin{equation}\label{eq:adjoint-local}
 w_t\ge (dw)_{qq}-(bw)_q
\end{equation}
with nonnegative parabolic trace.  Here this means that the initial trace is
nonnegative in $L^1(D)$ and the negative part has zero lateral trace in the
Dirichlet transposition sense.  If
$w(0,\xi)\ge\delta$ for almost every $\xi$ in an interval
$E\Subset D$ of positive length, then
$w(t,q)>0$ for every $t\in(0,t_*]$ and $q\in D$.
\end{lemma}
\begin{proof}
Choose $0\le z_0\le w(0,\cdot)$, supported in $E$ and equal to
$\delta/2$ on a smaller interval $E_0\Subset E$.  The bounded measurable
coefficients and uniform ellipticity give a unique nonnegative killed
solution in the transposition class,
\[
 z(t,q)=\int_DK_D(t,q;0,\xi)z_0(\xi)\,d\xi,
\]
by \cite[Theorem~4 and Corollary~4.1]{Aronson1968}.  The interior Harnack
principle \cite[Theorem~H]{Aronson1968} and the Gaussian bounds in
\cite{Aronson1967} give, for fixed $(t,q)$,
\begin{equation}\label{eq:killed-lower-general}
 \inf_{\xi\in E_0}K_D(t,q;0,\xi)>0;
\end{equation}
the constant depends only on the displayed structural bounds and the compact
sets.  Backward Dirichlet duality applied to the supersolution $w-z$ gives
$w\ge z$.  Consequently,
\[
 w(t,q)\ge z(t,q)
 \ge\frac\delta2\int_{E_0}K_D(t,q;0,\xi)\,d\xi>0
\]
by \eqref{eq:killed-lower-general}.
\end{proof}

\subsection{Proof of the strict-order theorem}
We give the complete argument because this is the step that rules out
disconnected active regions.

\emph{Step 1: derive the half-line equations from the whole-line inverse
equation.}
The unique weak
solution $\Xi$ of \eqref{eq:X-PDE} has the slope bounds
\eqref{eq:X-slope}.  On $p=q>0$ its flux factors as
\[
 \mathcal M(\tau,q,\Xi_p)=q\Phi^+(\tau,q\Xi_p)=q\Phi^+(\tau,F^+).
\]
Consequently, for every test function compactly supported in
$(0,T_0)\times(0,\infty)$,
\[
 F^+_\tau=q\partial_q^2\{q\Phi^+(\tau,F^+)\}
 =\partial_q\{q^2\partial_q\Phi^+(\tau,F^+)\}.
\]
The change $p=-q$ gives the same identity for $F^-$.  The initial traces are
those displayed below \eqref{eq:rho-barriers}, and
\eqref{eq:X-slope} gives the two linear-growth bounds in
\eqref{eq:F-PDE}.  Thus $F^\pm$ are restrictions of the unique whole-line
solution and are identified weak solutions of \eqref{eq:F-PDE}.

\emph{Step 2: nonstrict scaling order.}
Fix one sign and abbreviate $\Phi^\pm$ and $F^\pm$ by $\Phi$ and $F$.
Fix $c>1$ and set
 \[
 \begin{aligned}
 F_c(\tau,q)&=F(\tau,cq),\\
 d&=F_c-F,\\
 e&=\Phi(\tau,F_c)-\Phi(\tau,F)=A d.
 \end{aligned}
\]
Substitution into \eqref{eq:F-PDE} shows that $F_c$ solves the same
equation.  On every active branch,
$\partial_r\Phi(\tau,r)=a_i(\tau)/(2r^2)$; hence $F,F_c\asymp q$
and the secant formula give
\begin{equation}\label{eq:q2A}
 0<d_0\le q^2A(\tau,q)\le d_1<\infty.
\end{equation}
Applying \cref{lem:local-kato} gives
\begin{equation}\label{eq:kato-q}
 \partial_\tau d^-
 \le \partial_q\left[q^2\partial_q e^-\right],
 \qquad e^-=A d^- .
\end{equation}

The endpoints require an explicit check.  By joint continuity of $\Xi_p$
at $p=0$,
 \[
 \begin{aligned}
 \frac{d(\tau,q)}q={}&(c-1)\Xi_p(\tau,0)\\
 &+c[\Xi_p(\tau,cq)-\Xi_p(\tau,0)]\\
 &-[\Xi_p(\tau,q)-\Xi_p(\tau,0)].
 \end{aligned}
\]
The first term has a positive minimum on $[0,T_0]$ and the other terms
converge uniformly to zero.  Hence $d^-=0$ on $(0,q_\star)$ for some $q_\star>0$.
At infinity choose a $C^2$ weight $\omega$ satisfying
$\omega=1$ on $(0,1)$ and $\omega=q^{-4}$ for $q\ge2$.  Then
$\int\omega d^-<\infty$ and, by \eqref{eq:q2A},
\begin{equation}\label{eq:weight-q}
 \bigl(q^2\omega'(q)\bigr)'A(\tau,q)\le C\omega(q).
\end{equation}
Apply \eqref{eq:kato-q} first on $(\varepsilon,2R)$.  If
$\varepsilon<q_\star$, the inner boundary contribution is zero because
$d^-=e^-=0$ there.  Test with
$\omega_R(q)=\omega(q)\chi(q/R)$ and integrate twice by parts.  The terms
where $\chi'$ or $\chi''$ is nonzero are $o(1)$ as $R\to\infty$, because
$d^-=O(q)$, $e^-=O(q^{-1})$ on $R<q<2R$; the remaining term is controlled
by \eqref{eq:weight-q}.  Letting first $\varepsilon\downarrow0$ and then
$R\to\infty$ yields
\[
 \frac d{d\tau}\int_0^\infty\omega d^-\,dq
 \le C\int_0^\infty\omega d^-\,dq.
\]
The two initial derivatives are, without a sign convention left implicit,
 \[
 \begin{aligned}
 (F_0^+)'(q)&=\frac{R_G'(I_0(q))}{G''(I_0(q))}>0,\\
 (F_0^-)'(q)&=\frac{R_G'(I_0(-q))}{G''(I_0(-q))}>0.
 \end{aligned}
\]
Thus
$F_0(cq)>F_0(q)$, and Gronwall's inequality gives $F_c\ge F$.

\emph{Step 3: strictness.}
Put $W=F_c-F\ge0$ and define the secant coefficient $A_c$ by
$\Phi(F_c)-\Phi(F)=A_cW$.  Where $W=0$, choose any Borel value satisfying
the bounds in \eqref{eq:q2A}; the product $A_cW$ is independent of that
choice.  Then
\begin{equation}\label{eq:Wq}
 \begin{aligned}
 W_\tau&=\partial_q\left[q^2\partial_q(A_cW)\right]\\
 &=(D_cW)_{qq}-(b_cW)_q,\\
 D_c&=q^2A_c,\qquad b_c=2qA_c.
 \end{aligned}
\end{equation}
Equation \eqref{eq:q2A} makes $D_c$ uniformly positive and bounded; on
each $D=(q_-,q_+)\Subset(0,\infty)$, $b_c$ is bounded.  Fix
$(\tau_*,q_*)$ and choose $E\Subset D$ with $q_*\in D$.  The initial strict
order supplies $W(0,\xi)\ge\delta_E>0$ on $E$.  The function $W$ is
continuous and nonnegative, and its restriction to $D$ has nonnegative
lateral trace.  Applying \cref{lem:adjoint-positivity} directly to the
measurable coefficients $D_c,b_c$ gives $W(\tau_*,q_*)>0$.  Since the point
was arbitrary,
$F(\tau,cq)>F(\tau,q)$ for all $c>1$, $q>0$, and $\tau>0$.

\emph{Step 4: return to the state coordinate.}
On the positive branch,
$R(\tau,\Xi(\tau,q))=F^+(\tau,q)$.  On the negative branch,
$R(\tau,\Xi(\tau,-q))=-F^-(\tau,q)$.  Since $p=v_x$ and $\Xi$ are
increasing, the strict increase of $F^\pm$ proves strict increase of $R$
on both sides of its unique zero and hence on all of $\R$.  Finally,
\eqref{eq:second-derivative-bounds} and \eqref{eq:marginal-displacement} give
$R(\tau,x)\to\pm\infty$ as $x\to\pm\infty$.

\begin{remark}
The theorem proves strict level crossing, not transversality.  A continuous
strictly increasing function may have zero derivative at its unique crossing.
Accordingly, the main active-set theorem asserts only continuity of switching
locations; a separate remark below records sufficient conditions for $C^1$
regularity.
\end{remark}
}

\section{Construction and verification of the ordered feedback}\label{sec:hull}
For fixed $\tau$, associate action $i$ with the generator point
\[
 P_i(\tau)=\left(\beta_i(\tau),\frac{a_i(\tau)}2\right)
\]
and the line
\begin{equation}\label{eq:line}
 L_i(\tau,r)=\frac{a_i(\tau)}2-\beta_i(\tau)r.
\end{equation}
Because $v_{xx}>0$, action selection in \eqref{eq:hjb} is equivalent to minimizing $L_i(\tau,R(\tau,x))$.

Before taking a hull, identical and dominated labels must be removed.  Actions with identical $(\beta_i,a_i)$ are equivalent and may be represented by one label.  Among actions with the same $\beta_i$, only the one with smallest $a_i$ can be active.  After this preprocessing, let
\[
 P_{j_1}(\tau),\ldots,P_{j_K}(\tau),
 \qquad \beta_{j_1}<\cdots<\beta_{j_K},
\]
be the vertices of the lower convex hull, with collinear interior points omitted.
Set $\Delta\beta_k=\beta_{j_{k+1}}-\beta_{j_k}$.  Its adjacent slopes are
\begin{equation}\label{eq:kappa}
 \kappa_k(\tau)=
 \frac{a_{j_{k+1}}(\tau)-a_{j_k}(\tau)}
 {2\Delta\beta_k(\tau)},
 \qquad k=1,\ldots,K-1.
\end{equation}
Convexity of the hull gives $\kappa_1<\cdots<\kappa_{K-1}$.  A nonvertex strictly above the hull is never active; a discarded collinear point can tie only at one threshold and has no open active region.

\begin{theorem}[Ordered active set]\label{thm:active-set}
Under \cref{ass:main}, fix $\tau>0$.  For each $k$ there is exactly one
$s_k(\tau)\in\R$ satisfying
\begin{equation}\label{eq:switch-level}
 R(\tau,s_k(\tau))=\kappa_k(\tau).
\end{equation}
The switching points are ordered,
$s_1(\tau)<\cdots<s_{K-1}(\tau)$, and the minimizing hull vertex is
\begin{equation}\label{eq:ordered-policy}
 j_1,\ j_2,\ldots,j_K
\end{equation}
on the successive state intervals.  The maps $\tau\mapsto s_k(\tau)$ are
continuous as long as the hull combinatorics do not change.
\end{theorem}
\begin{proof}
Existence, uniqueness, and ordering follow from \cref{thm:strict-order} and
the ordered hull slopes.  The lower-envelope property of convex hulls gives
\eqref{eq:ordered-policy}.  Continuity follows from locally uniform
continuity of $R$ and strict crossing.
\end{proof}

For an operations manager, the theorem converts a potentially complicated
state-feedback table into a sorted threshold policy.  Offline computation
stores $K-1$ curves.  Online implementation compares the observed state with
those curves and selects the corresponding hull action; it never evaluates the
HJB equation in real time.  At a given time only that time's hull actions need
to be represented; a label may be removed from the full-horizon table only
when it is absent from every pointwise hull.  For time-invariant actions this
reduces to one exact offline hull calculation.  This compression differs
from approximate aggregation of nearby control levels.

The H\"older assumptions of the main theorem imply continuity of the switching
curves, but not differentiability: strict level crossing alone does not imply
transversality.

The preceding theorem identifies a pointwise minimizer.  The next result supplies the missing stochastic verification.
\begin{theorem}[Weak Markov verification]\label{thm:verification}
Let $i^*(\tau,x)$ be the smallest label attaining the minimum in
\eqref{eq:hjb}; any other fixed Borel tie-breaking rule gives the same conclusion.  Then:
\begin{enumerate}
\item $i^*$ is Borel measurable;
\item for every $(\tau,x)$ there exists a weak solution of
\begin{equation}\label{eq:closed-loop}
 \begin{aligned}
 dX_s={}&-\beta_{i^*(\tau-s,X_s)}(\tau-s)\,ds\\
 &+\sqrt{a_{i^*(\tau-s,X_s)}(\tau-s)}\,dW_s.
 \end{aligned}
\end{equation}
\item every weak solution of \eqref{eq:closed-loop} is optimal for
\eqref{eq:value}.  In particular, an optimal weak Markov control exists.
\end{enumerate}
\end{theorem}
\begin{proof}
For $i\in\mathcal I$, set
\[
 g_i(\tau,x)=-\beta_i(\tau)v_x(\tau,x)
 +\frac{a_i(\tau)}2v_{xx}(\tau,x).
\]
By \cref{cor:regularity}, every $g_i$ is continuous on
$(0,T]\times\R$.  The smallest-index minimizer of finitely many continuous
functions is Borel; define it arbitrarily on the terminal slice.  The
closed-loop drift and variance are therefore bounded Borel functions, and
the variance lies in $[\underline a,\overline a]$.  Standard weak existence
for uniformly nondegenerate diffusions with bounded measurable coefficients
gives a solution of \eqref{eq:closed-loop}; see
\cite{StroockVaradhan1979,Krylov1980}.

It remains to verify optimality without assuming strong uniqueness.  Fix an
arbitrary admissible weak system, write
$u(s,z)=v(\tau-s,z)$, and, for $n\ge1$ and $0<\delta<\tau$, stop at
\[
 \sigma_{n,\delta}=(\tau-\delta)
 \wedge\inf\{s\ge0:|X_s|\ge n\}.
\]
On the corresponding positive-time compact cylinder, $u$ is $C^{1,2}$ and
the Bellman equation gives, for every action $i$,
\[
 u_s-\beta_i(\tau-s)u_x
 +\frac{a_i(\tau-s)}2u_{xx}\ge0.
\]
It\^o's formula and bounded stopping thus yield
\[
 \E\bigl[v(\tau-\sigma_{n,\delta},X_{\sigma_{n,\delta}})\bigr]
 -v(\tau,x)
 =\E\!\int_0^{\sigma_{n,\delta}}
 \left(u_s-\beta_{I_s}u_x+\frac{a_{I_s}}2u_{xx}\right)(s,X_s)\,ds
 \ge0.
\]
Bounded coefficients give a fourth-moment estimate for
$\sup_{s\le\tau}|X_s|$.  Together with the quadratic growth of $v$, this
provides uniform integrability.  Letting first $n\to\infty$ and then
$\delta\downarrow0$, using the locally uniform terminal convergence
$v(\delta,\cdot)\to G$ and the same growth bound, gives
$\E[G(X_\tau)]\ge v(\tau,x)$.

For a weak solution of \eqref{eq:closed-loop}, the displayed generator
inequality is an equality almost everywhere in time.  The same localization
limits therefore give $\E[G(X_\tau)]=v(\tau,x)$.  This argument depends only
on the selected generator, so any fixed Borel rule on ties is optimal.
\end{proof}

\section{Operational action screening and continuous benchmarks}\label{sec:scope}

\subsection{Exact screening of a finite action set}
Let $\mathcal H(\tau)$ be the ordered set of vertices of the lower convex
hull of
\[
 \mathcal P(\tau)=\{(\beta_i(\tau),a_i(\tau)/2):i\in\mathcal I\}.
\]
Theorem~\ref{thm:active-set} first gives a pointwise-in-time screening rule:
at time $\tau$, every action outside $\mathcal H(\tau)$ is suboptimal for all
states.  Define the all-time effective action set
\[
 \mathcal I_{\rm eff}:=\bigcup_{0\le\tau\le T}\mathcal H(\tau),
\]
where the union is over action labels.  Only labels outside
$\mathcal I_{\rm eff}$ may be removed permanently before solving the dynamic
program.  If the action characteristics are time independent, the hull is
fixed and the single hull calculation is an exact preprocessing step.  This
distinction matters: an action can lie above today's hull and become a hull
vertex later.

The case of uncontrolled diffusion is contained in the state model by taking
$a_i\equiv a$, but its policy proof should not be described as a mechanical
application of the general ratio theorem.  The common second-order term
cancels from every action comparison, so the decision depends directly on the
sign of $v_x$.  Strict convexity, rather than strict monotonicity of
$v_x/v_{xx}$, supplies the ordering.

\begin{corollary}[Drift-only and variance-only action sets]\label{cor:special-menus}
Under \cref{ass:main}, the following conclusions hold.
\begin{enumerate}
\item If $a_i(\tau)\equiv a(\tau)$ for every action, then all intermediate
drift levels are redundant.  Where $v_x<0$ the smallest $\beta_i$ is optimal,
and where $v_x>0$ the largest $\beta_i$ is optimal.  Since $v_x$ is strictly
increasing and onto, these two regions are separated by exactly one state
$m(\tau)$ determined by $v_x(\tau,m(\tau))=0$.
\item If $\beta_i(\tau)\equiv\beta(\tau)$ for every action, an action with minimum
$a_i(\tau)$ is optimal at every state; all larger-variance actions are
redundant.
\item If both coordinates vary, precisely the lower-hull vertices may occupy
open optimality regions, in their hull order.
\end{enumerate}
\end{corollary}
\begin{proof}
If $a_i\equiv a$, then for actions $i,j$,
\[
 \left(-\beta_i v_x+\frac a2v_{xx}\right)
 -\left(-\beta_j v_x+\frac a2v_{xx}\right)
 =-(\beta_i-\beta_j)v_x.
\]
Proposition~\ref{prop:second-derivative} gives $v_{xx}>0$ for positive time, so
$v_x(\tau,\cdot)$ is strictly increasing.  Lemma~\ref{lem:displacement}
and the fact that $G'$ is onto imply that $v_x(\tau,\cdot)$ is onto as well.
It therefore has one zero $m(\tau)$.  For $x<m(\tau)$ the smallest drift
coefficient minimizes the Hamiltonian; for $x>m(\tau)$ the largest one does.
No ratio equation is needed.  If instead $\beta_i\equiv\beta$, then
\[
 \left(-\beta v_x+\frac {a_i}2v_{xx}\right)
 -\left(-\beta v_x+\frac {a_j}2v_{xx}\right)
 =\frac{a_i-a_j}{2}v_{xx},
\]
and positivity of $v_{xx}$ selects minimum variance at every state.  The
jointly controlled case is Theorem~\ref{thm:active-set}.
\end{proof}

For operations design, the corollary distinguishes a large nominal action set
from its effective subset.  An action above the lower hull is never uniquely optimal
under the model and can be removed from policy storage and online
implementation.  If the lower hull contains $m$ vertices, the feedback is
encoded by $m-1$ state thresholds rather than by a general state-to-action map.
This conclusion concerns operating decisions for a given action set.

\subsection{A zero-gap relaxation and a genuinely different continuum}
The word ``continuous approximation'' can refer to two different models.
The first is a relaxed generator control.  At each time it chooses
\[
 (\bar\beta,\bar a)=\sum_i\lambda_i(\beta_i,a_i),
 \qquad \lambda_i\ge0,\qquad \sum_i\lambda_i=1.
\]

\begin{proposition}[Zero gap for generator convexification]\label{prop:zero-gap}
Replacing $\mathcal P(\tau)$ by $\operatorname{co}\mathcal P(\tau)$ leaves
the Hamiltonian and the value function unchanged.
\end{proposition}
\begin{proof}
The Hamiltonian is linear in the mean--variance point, so its minimum over a
convex hull is attained at an original action point.  The HJB equations and
terminal data coincide, and comparison gives equality of the values;
equivalently, relaxed controls and chattering give the same conclusion
\cite{ElKarouiEtAl1987,HaussmannLepeltier1990}.
\end{proof}

This relaxation averages variance rates, not diffusion amplitudes.  In
general,
\[
 \left(\sum_i\lambda_i\sqrt{a_i}\right)^2
 \ne \sum_i\lambda_i a_i.
\]
Consequently, interpolating actuator amplitudes defines a different physical
model and is not covered by \cref{prop:zero-gap}.

To make an implementable comparison, suppose a physical continuum
$U=[\underline u,\overline u]$ is specified with drift and variance
$\beta(\tau,u)$ and $a(\tau,u)$, and the available finite action set is
$U_\Delta=\{u_1,\ldots,u_N\}\subset U$.  For a progressively measurable
 $U$-valued process $u$, define
 \[
 \begin{aligned}
 X_s^{u}={}&x-\int_0^s\beta(\tau-r,u_r)\,dr\\
 &+\int_0^s\sqrt{a(\tau-r,u_r)}\,dW_r,\\
 v^c(\tau,x)&=\inf_u\E G(X_\tau^u).
 \end{aligned}
\]
The value $v^\Delta$ is defined by restricting $u$ to $U_\Delta$.  If a
Borel feedback $\widehat u$ has a weak closed-loop solution, write
$J^{\widehat u}(\tau,x)=\E G(X_\tau^{\widehat u})$.  If
$u^c(\tau,x)$ is a measurable continuous-model optimizer, its nearest-action
implementation is
\begin{equation}\label{eq:nearest-quantization}
 \widehat u^{\rm near}(\tau,x)
 \in\arg\min_{u_i\in U_\Delta}|u_i-u^c(\tau,x)|.
\end{equation}
It need not equal the exact finite policy, which is computed from
$v^\Delta$.  Let $J^{\rm near}=J^{\widehat u^{\rm near}}$.  Admissibility
and set inclusion give
\begin{equation}\label{eq:value-benchmark-chain}
 v^c\le v^\Delta\le J^{\rm near}.
\end{equation}
Thus the relevant implementation loss is $J^{\rm near}-v^\Delta$, not the
zero gap between the finite problem and its generator convexification.

The next result makes the effect of action-set resolution explicit.  Let
\begin{equation}\label{eq:menu-radius}
 \Delta:=\sup_{u\in[\underline u,\overline u]}
          \min_{u_i\in U_\Delta}|u-u_i|
\end{equation}
be the covering radius of the finite action set, and assume that both endpoints of the
continuous interval belong to $U_\Delta$.

\begin{proposition}[Quadratic action-set approximation bound]
\label{prop:menu-error}
Suppose $\beta,a\in C^{\alpha/2}([0,T];C^2(U))$, their derivatives in $u$
through order two are bounded, and $a\ge\underline a>0$.  Assume that the
continuous problem has a classical value satisfying
\[
 |v_x^c(\tau,x)|\le C(1+|x|),
 \qquad |v_{xx}^c(\tau,x)|\le C,
\]
and a Borel minimizer $u^c$.  Let $Q_\Delta$ be a Borel nearest-action map,
assume that both endpoints of $U$ belong to $U_\Delta$, and assume that the
closed-loop martingale problem under
$\widehat u^\Delta=Q_\Delta(u^c)$ has a weak solution.  Set
\[
 h^c(\tau,x,u)
 =-\beta(\tau,u)v_x^c(\tau,x)
   +\frac{a(\tau,u)}2v_{xx}^c(\tau,x).
\]
Then, for $0\le\tau\le T$,
\begin{equation}\label{eq:quadratic-menu-bound}
 \begin{aligned}
 0&\le v^\Delta(\tau,x)-v^c(\tau,x)\\
 &\le J^{\rm near}(\tau,x)-v^c(\tau,x)\\
 &\le C_T(1+|x|)\Delta^2,
 \end{aligned}
\end{equation}
and consequently
\begin{equation}\label{eq:quadratic-policy-loss}
 0\le J^{\rm near}(\tau,x)-v^\Delta(\tau,x)
 \le C_T(1+|x|)\Delta^2.
\end{equation}
The constant depends on the coefficient and value-derivative bounds and on
$T$, but not on the finite action set.
\end{proposition}
\begin{proof}
The assumptions give the checkable bound
\begin{equation}\label{eq:hessian-u-bound}
 \begin{aligned}
 h^c_{uu}&=-\beta_{uu}v_x^c
 +\frac12a_{uu}v_{xx}^c,\\
 |h^c_{uu}(\tau,x,u)|&\le K(1+|x|).
 \end{aligned}
\end{equation}
Fix $(\tau,x)$.  If $u^c$ is an interior minimizer, then
$h_u^c(\tau,x,u^c)=0$.  If it is an endpoint minimizer, the endpoint belongs
to the finite action set and $Q_\Delta(u^c)=u^c$.  Taylor's formula and
\eqref{eq:hessian-u-bound} therefore give in both cases
\begin{equation}\label{eq:local-menu-gap}
 \begin{aligned}
 0&\le h^c(\tau,x,Q_\Delta(u^c))
 -h^c(\tau,x,u^c)\\
 &\le\frac K2(1+|x|)\Delta^2.
 \end{aligned}
\end{equation}
Let $X^{\Delta}$ be a weak solution under the rounded Markov policy.
Because the coefficients are bounded on the compact control interval and do
not grow with the state,
\[
 \E\left[\sup_{0\le s\le\tau}|X_s^{\Delta}|\right]
 \le C_T(1+|x|).
\]
Apply It\^o's formula to
$v^c(\tau-s,X_s^{\Delta})$, localize, and use the quadratic growth of the
value function to remove the localization.  Since
$v_\tau^c=\min_u h^c(\tau,x,u)$, we obtain
\begin{align*}
 J^{\rm near}-v^c
 &=\E\int_0^\tau
 \bigl[h^c(\tau-s,X_s^{\Delta},\widehat u_s^\Delta)\\
 &\hspace{45pt}-\min_u h^c(\tau-s,X_s^{\Delta},u)\bigr]ds\\
 &\le C_T(1+|x|)\Delta^2
\end{align*}
by \eqref{eq:local-menu-gap}.  Finally, the continuous action set contains
the finite action set, while the rounded policy is finite-action admissible.  Hence
 $v^c\le v^\Delta\le J^{\rm near}$, which proves both displayed bounds.
\end{proof}

The quadratic rate uses first-order optimality of the smooth continuous
Hamiltonian.  With only Lipschitz dependence of drift and diffusion on the
control, direct coupling generally gives only an $O(\Delta)$ bound.  In the
production model \eqref{eq:operational-continuum-new},
$h^c_{uu}=v_{xx}^c$; hence the derivative bound assumed in
\cref{prop:menu-error} verifies \eqref{eq:hessian-u-bound} without an
additional state-growth factor.  The finite-action second-derivative theorem is not
being invoked for the continuous action set at this point.

\section{Numerical experiments and operational implications}\label{sec:numerics}
The experiments serve three purposes.  The first quantifies the cost
of using a continuously optimized feedback in a system that can implement
only finitely many actions.  The second illustrates exact lower-hull screening and
ordered regions; the third tests pointwise screening with time-varying actions
and a strongly asymmetric loss.  All calculations use backward Euler time
stepping and Howard policy iteration.  For each experiment, we report grid,
domain, residual, and policy-iteration checks.

\subsection{Continuous planning versus finite implementation}
Consider the physical actuator
\begin{equation}\label{eq:operational-continuum-new}
 \begin{gathered}
 dX_s=-u_s\,ds+\sqrt{0.4+u_s^2}\,dW_s,\\
 \frac12\le u_s\le2,
 \qquad G(x)=\frac{x^2}{2},
 \end{gathered}
\end{equation}
over a unit horizon.  The additive variance $0.4$ represents background
uncertainty that remains under every operating decision.  The continuous
model's feedback is
\begin{equation}\label{eq:continuous-feedback-new}
 u^c(\tau,x)=\Pi_{[1/2,2]}
 \left(\frac{v_x^c(\tau,x)}{v_{xx}^c(\tau,x)}\right).
\end{equation}
An equally spaced action set with two, three, or five entries is then the actual
implementable action set.  This is the dynamic analogue of receiving a
continuous recommendation of 5.5 units when only five or six can be used:
nearest rounding is feasible, but its threshold is computed from a different
continuation value than the exact finite-action threshold.
The positive lower rate $u=1/2$ represents a batch or service operation that
must remain online at a certified minimum rate once the planning window has
started; an idle action is intentionally not available in this benchmark.

\begin{figure}[t]
\centering
\includegraphics[width=0.98\textwidth]{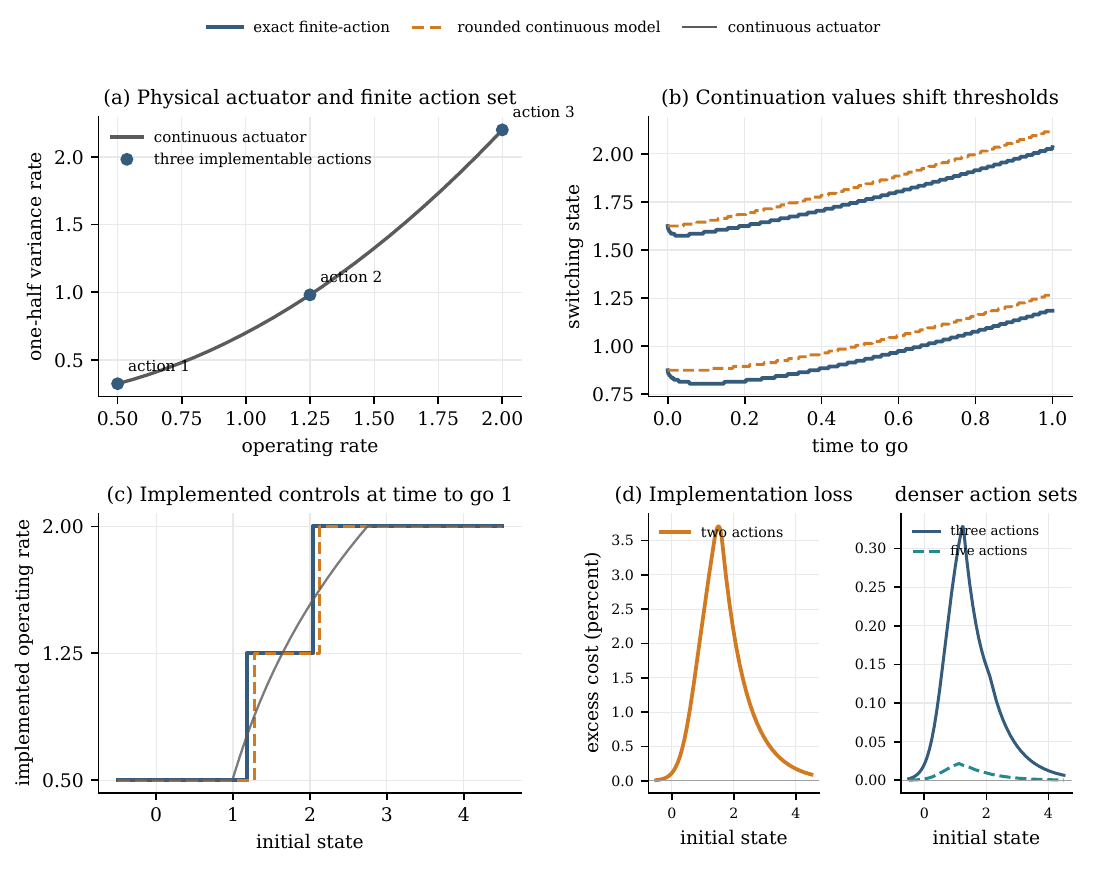}
\caption{Continuous planning and finite implementation for
\eqref{eq:operational-continuum-new}.  (a) The physical continuous actuator
and a three-action set.  (b) Switching states of the exact finite-action policy
and of the policy that selects the nearest feasible action to the continuous
feedback.  (c) The three feedbacks at one unit of time to go.  (d) Relative
implementation loss; the two adjacent axes use labeled scales so that the
coarse two-action loss and the smaller three- and five-action losses remain fully
visible without an overlay.}
\label{fig:operational-comparison}
\end{figure}

On the grid $\Delta x=10^{-2}$, $\Delta\tau=5\times10^{-4}$ over
$[-6,8]$, the peak relative loss
$100(J^{\rm near}-v^\Delta)/v^\Delta$ on $[-1,5]$ is $3.705\%$ for the
coarse two-action set.  Refining to three and five equally spaced actions reduces
the peak to $0.329\%$ and $0.0217\%$, respectively.  The first number
identifies sparse action sets as the regime in which exact finite-action optimization
matters most; the latter two are consistent with the squared-radius upper
bound and illustrate how rapidly implementation risk declines as the action
set is refined.

A coarse sensitivity audit clarifies the mechanism.  Reducing background
variance from $0.4$ to $0.2$ raises the two-action peak loss from $3.67\%$ to
$5.72\%$, whereas increasing it to $0.8$ lowers the peak to $1.82\%$.
Increasing the rate--variance coupling from $0.6$ to $1.4$ raises the peak
absolute loss from $0.00825$ to $0.04322$.  Extending the horizon from $0.5$
to $1.5$ raises the peak relative loss from $2.45\%$ to $4.23\%$.  These are
comparative numerical findings, not general comparative statics: changing a
parameter changes both the continuous feedback and the finite continuation
value.

\subsection{Exact screening and ordered regions}
For a nonquadratic illustration, take
\begin{equation}\label{eq:numerical-HJB-new}
 \begin{aligned}
 v_\tau&=\min_{1\le i\le4}
 \left\{-\beta_i v_x+\frac{a_i}{2}v_{xx}\right\},\\
 G(x)&=\frac{x^2}{2}
 +0.6\bigl(\sqrt{1+x^2}-1\bigr),
 \end{aligned}
\end{equation}
with
$(\beta_i,a_i)=(1,1.25),(1.5,2.5),(1.75,4.5),(2,4.25)$.
The third generator point lies strictly above the segment joining the second
and fourth points.  It is therefore screened before solving the HJB equation,
and the two remaining hull slopes are $1.25$ and $1.75$.  The computed
switching states at one unit of time to go are $2.021962$ and $2.498058$; enlarging
the domain changes them by less than $1.3\times10^{-5}$.

\begin{figure}[!ht]
\centering
\includegraphics[width=0.98\textwidth]{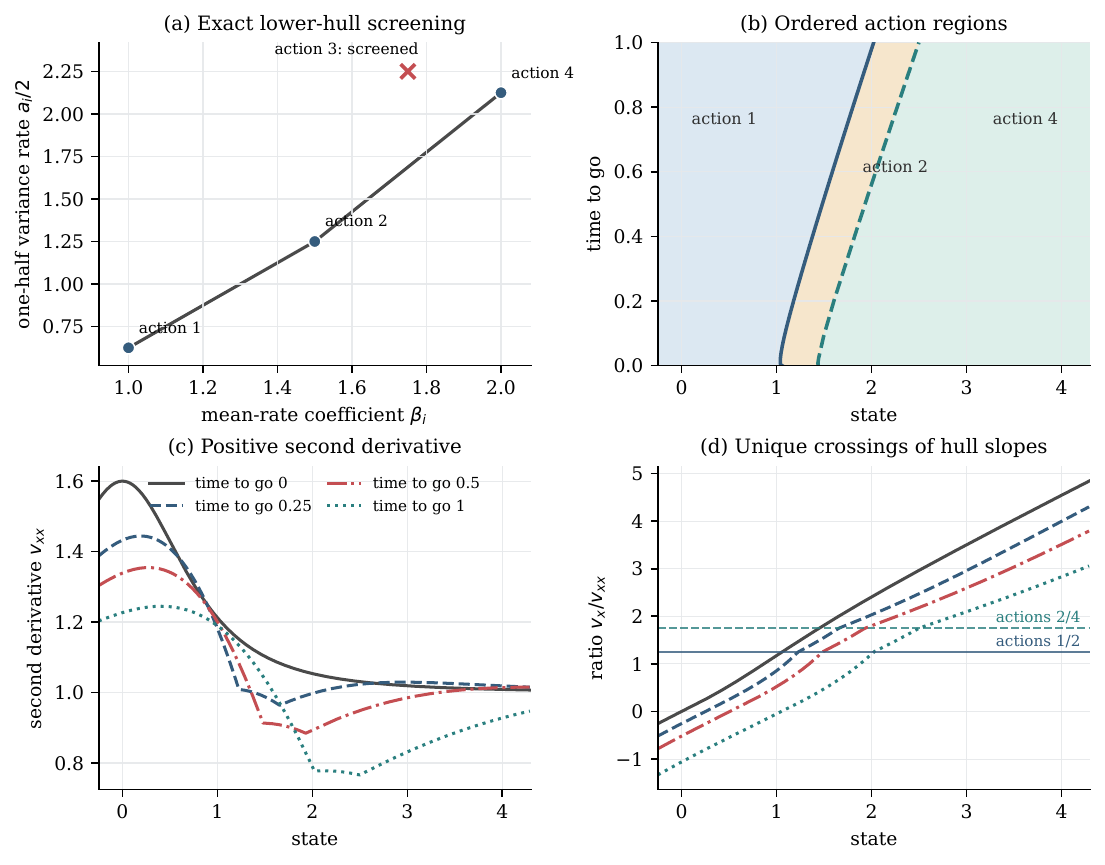}
\caption{Four-action benchmark \eqref{eq:numerical-HJB-new}.  (a) Action 3 is
strictly above the lower hull.  (b) The computed regions follow the hull
order and never use action 3.  (c) The computed second derivative remains
positive at the displayed times.  (d) The increasing ratio $v_x/v_{xx}$
crosses each lower-hull slope exactly once.}
\label{fig:numerics}
\end{figure}

The screen is also inexpensive for a large action set.  In a synthetic
64-action audit, 47 actions are removed, leaving 17 hull actions and 16 stored thresholds.
Sorting and the hull scan require $O(N\log N)$ time, or $O(N)$ for an already
sorted set, while online binary search requires $O(\log K)$ comparisons for
$K$ retained actions.  Timing details and the full audit accompany the code.

\subsection{Time-varying screening and an asymmetric terminal loss}
The pointwise screen matters when action characteristics change with time.  On
a three-unit horizon we retain the four drift rates in
\eqref{eq:numerical-HJB-new}, keep the variances of actions 1, 2, and 4 fixed,
and set
\[
 \frac{a_3(\tau)}2=1.665+0.085\cos(2\pi\tau/3).
\]
The height of action 3 above the chord joining actions 2 and 4 is therefore
$-0.0225+0.085\cos(2\pi\tau/3)$: action 3 is screened near both ends of the
horizon and becomes a lower-hull vertex in the middle.

To model terminal shortfall as substantially more costly than surplus, take
\begin{equation}\label{eq:asymmetric-terminal}
 \begin{aligned}
 G_\delta(x)=\frac{x^2}{2}+\frac\delta2\bigl\{
 x\sqrt{1+x^2}+\operatorname{arsinh}x-2x\bigr\},\\
 \delta=0.8.
 \end{aligned}
\end{equation}
Its left and right asymptotic quadratic weights differ by a factor of nine.
Direct calculation gives
\[
 \begin{aligned}
 G_\delta'(x)&=x+\delta(\sqrt{1+x^2}-1),\\
 G_\delta''(x)&=1+\frac{\delta x}{\sqrt{1+x^2}},\qquad
 G_\delta'''(x)=\frac\delta{(1+x^2)^{3/2}}.
 \end{aligned}
\]
Hence $1-\delta\le G_\delta''\le1+\delta$.  For $x\le0$ one has
$G_\delta'(x)\le0$ and thus $R_{G_\delta}'(x)\ge1$.  For $x\ge0$,
$G_\delta'(x)\le(1+\delta)x$ and $G_\delta''(x)\ge1$, so
\[
 R_{G_\delta}'
 \ge1-\frac{\delta(1+\delta)x}{(1+x^2)^{3/2}}
 \ge1-\frac{2\delta(1+\delta)}{3\sqrt3}>0.445.
\]
Here the last maximum is attained at $x=1/\sqrt2$.  Since
$G_\delta''\ge1-\delta$ and the displayed derivatives are bounded in the
required way, $R_{G_\delta}'$ also has a finite global upper bound.  Thus
\cref{ass:main} holds, with zero as the unique minimizer.

\begin{figure}[t]
\centering
\includegraphics[width=0.98\textwidth]{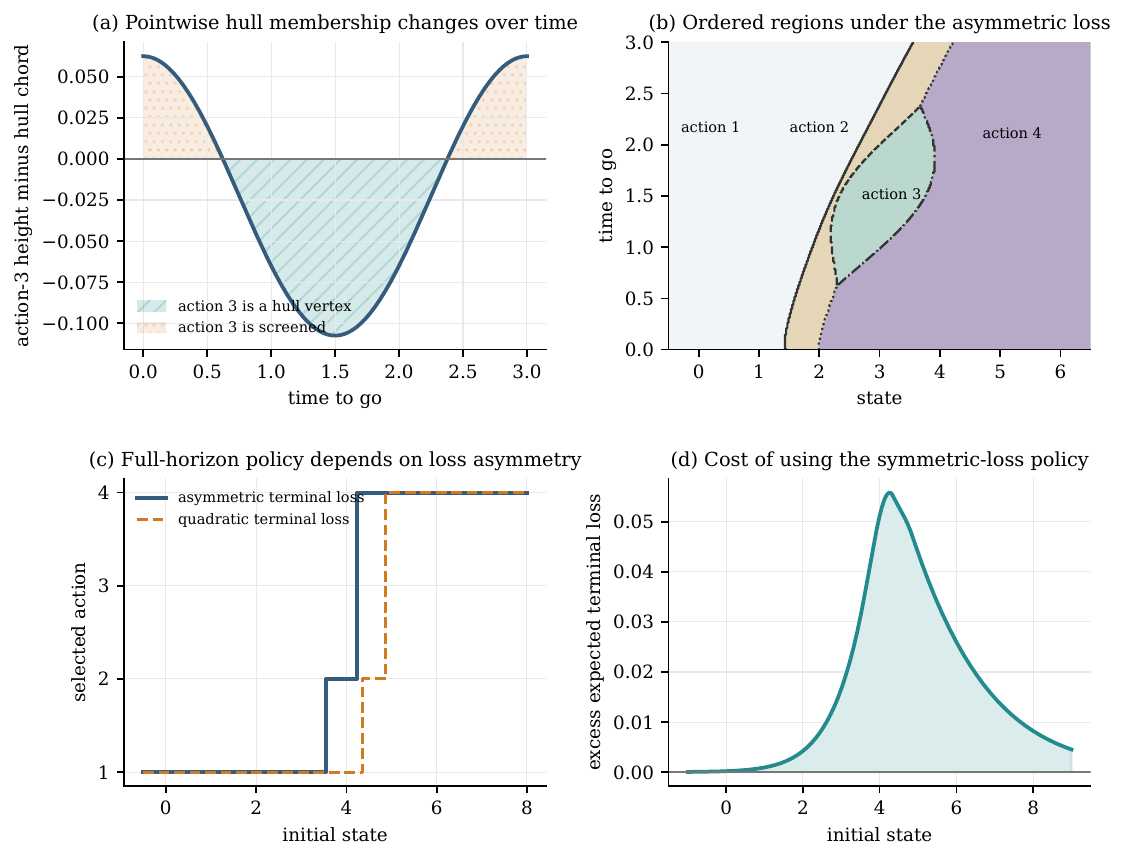}
\caption{Time-varying screening and asymmetric terminal regulation.  (a) The
signed height of action 3 above the chord joining actions 2 and 4 determines its
pointwise hull membership.  (b) The computed optimal regions use action 3
exactly during its hull interval, apart from the two collinearity instants.
(c) The asymmetric loss \eqref{eq:asymmetric-terminal} shifts the full-horizon
thresholds relative to the quadratic loss.  (d) Applying the policy optimized
for the quadratic loss when \eqref{eq:asymmetric-terminal} is the true loss
produces a positive excess terminal cost.}
\label{fig:timevarying-asymmetric}
\end{figure}

On $[-12,18]$ with $\Delta x=0.02$ and $\Delta\tau=10^{-3}$, the selector has
no pointwise-screening mismatch away from a $3\times10^{-3}$ neighborhood of
the two collinearity instants.  The minimum computed $v_{xx}$ on
$[-1,9]$ is $0.272$, and the maximum Bellman residual is
$1.21\times10^{-9}$.  Using the policy optimized for the symmetric quadratic
loss raises expected terminal cost by $0.05582$ at its peak; normalized by
$1+|V^\star|$, the peak loss is $1.713\%$.  The coarser grid
$(0.04,0.002)$ gives $0.05596$ and $1.716\%$, respectively.  Thus the third
experiment simultaneously verifies pointwise screening and shows that a
policy calibrated to the wrong terminal asymmetry can be materially
suboptimal.

\section{Conclusion}
For a finite action set, the lower generator hull is an exact screen and
the surviving actions form an ordered threshold policy.  This converts the
dynamic program into a compact implementation table and applies pointwise to
time-varying action sets.  A convexified generator has zero relaxation gap, but a
different continuous actuator and nearest-action implementation need not be
optimal for the available action set; the quadratic action-set approximation bound and the
experiments quantify this distinction.  Drift-only and diffusion-only action sets
follow as simpler one-coordinate cases, and a common affine state drift is
removable.  Action-dependent state multipliers, within-horizon operating or
changeover costs, and higher-dimensional states require different structural
arguments and remain natural directions for action-set design and control.

\FloatBarrier

\appendix
\section{Second-derivative bounds: approximation and distributional passage}\label{app:second-derivative}
This appendix completes the proof of \cref{prop:second-derivative} and makes explicit how the distributional estimate becomes a pointwise one.

\subsection{Smooth Bellman approximation}
Let
\[
 \begin{aligned}
 \mathcal F^\varepsilon(\tau,p,r)
 =-\varepsilon\log\sum_{i=1}^N
 \exp\biggl(&\frac{\beta_i^\varepsilon(\tau)p}{\varepsilon}\\
 &-\frac{a_i^\varepsilon(\tau)r}{2\varepsilon}\biggr),
 \end{aligned}
\]
where the time coefficients are smoothly approximated without changing the ellipticity bounds.  Then
\begin{equation}\label{eq:Fe-bounds}
 \frac{\underline a}{2}\le \mathcal F^\varepsilon_r
 \le\frac{\overline a}{2},
 \qquad |\mathcal F^\varepsilon_p|\le\overline\beta,
\end{equation}
and $\mathcal F^\varepsilon\to\mathcal F$ locally uniformly.  Smooth $G^\varepsilon$ can be chosen with
$m\le G^{\varepsilon\prime\prime}\le M$.

For each fixed $\varepsilon$, standard whole-space Cauchy theory for smooth
uniformly parabolic Bellman equations with quadratic growth gives a classical
solution; see \cite{Krylov1980,LadyzhenskayaEtAl1968}.  Equivalently, one may
use finite-cylinder exhaustion with restrictions of common quadratic
sub- and supersolutions as lateral data; the interior estimates are
independent of the cylinder radius.  This exhaustion defines the lateral
data explicitly and permits passage to the whole-space solution.
Differentiating twice and writing
$w=v^\varepsilon_{xx}$ gives
\begin{align*}
 w_\tau={}&\mathcal F^\varepsilon_r w_{xx}
 +\mathcal F^\varepsilon_p w_x\\
 &+\mathcal F^\varepsilon_{pp}w^2
 +2\mathcal F^\varepsilon_{pr}ww_x\\
 &+\mathcal F^\varepsilon_{rr}w_x^2.
\end{align*}
At an interior first zero of $w$, one has $w_x=0$; the drift and nonlinear
terms vanish, while the diffusion term is nonnegative.  On the whole line
apply the argument to
$w+\delta e^{C_\varepsilon\tau}(1+x^2)$.  Polynomial-growth estimates allow
$C_\varepsilon$ to be chosen so that a negative first contact cannot escape
to spatial infinity.  Letting $\delta\downarrow0$ gives $w\ge0$ for the fixed
smooth approximation.  Dependence of $C_\varepsilon$ on the smoothing is
harmless: the uniform quantitative estimates below come from the equation
for $v_x$ rather than this barrier argument.

Viscosity stability gives $v^\varepsilon\to v$ locally uniformly.  For every nonnegative test function $\varphi$,
 \[
 \begin{aligned}
 \langle D_{xx}v,\varphi\rangle
 &=\int v\varphi_{xx}\\
 &=\lim_{\varepsilon\downarrow0}
 \int v^\varepsilon\varphi_{xx}\\
 &=\lim_{\varepsilon\downarrow0}
 \int v^\varepsilon_{xx}\varphi\ge0.
 \end{aligned}
\]
This is precisely the meaning of $D_{xx}v\ge0$ in distributions.  In one dimension it is equivalent to convexity of a continuous representative of $v(\tau,\cdot)$.

\subsection{Quantitative strictness}
For the smooth approximation, define $p^\varepsilon=v_x^\varepsilon$ and
$z_h^\varepsilon=p^\varepsilon(\tau,x+h)-p^\varepsilon(\tau,x)$.
Since
$p^\varepsilon_\tau=\partial_x\mathcal F^\varepsilon
(\tau,p^\varepsilon,p_x^\varepsilon)$, the two-variable mean-value formula
gives
\begin{equation}\label{eq:zh}
 \begin{aligned}
 (z_h^\varepsilon)_\tau
 &=\partial_x\bigl(A_h^\varepsilon(z_h^\varepsilon)_x
 +B_h^\varepsilon z_h^\varepsilon\bigr),\\
 \frac{\underline a}{2}\le A_h^\varepsilon
 &\le\frac{\overline a}{2},
 \qquad |B_h^\varepsilon|\le\overline\beta .
 \end{aligned}
\end{equation}
For example, if $\Delta_hf=f(\tau,x+h)-f(\tau,x)$, one may take
\begin{align*}
 A_h^\varepsilon
 &=\int_0^1\mathcal F_r^\varepsilon
 \bigl(\tau,p^\varepsilon(x+h),\\
 &\hspace{55pt}p_x^\varepsilon(x)
 +\theta\Delta_hp_x^\varepsilon\bigr)\,d\theta,\\
 B_h^\varepsilon
 &=\int_0^1\mathcal F_p^\varepsilon
 \bigl(\tau,p^\varepsilon(x)
 +\theta\Delta_hp^\varepsilon,\\
 &\hspace{55pt}p_x^\varepsilon(x)\bigr)\,d\theta,
\end{align*}
after splitting the flux difference first in its second and then in its
first argument.  The fundamental solution $\Gamma_h^\varepsilon$ obeys,
for $0<\tau\le T_0$,
 \[
 \begin{aligned}
 c\tau^{-1/2}e^{-|x-y|^2/(c\tau)}
 &\le\Gamma_h^\varepsilon(\tau,x;0,y)\\
 &\le C\tau^{-1/2}e^{-|x-y|^2/(C\tau)}.
 \end{aligned}
\]
Although the lower-order term in \eqref{eq:zh} is only bounded measurable, these estimates remain valid \cite{Aronson1967,Aronson1968}.  Integrating in $y$ and using
$mh\le z_h^\varepsilon(0,y)\le Mh$ yields constants independent of $\varepsilon,h,x$:
\[
 c_{T_0}mh\le z_h^\varepsilon(\tau,x)
 \le C_{T_0}Mh.
\]
Pass first $\varepsilon\downarrow0$ in distributions, then $h\downarrow0$.  We obtain
$c_{T_0}m\le D_{xx}v\le C_{T_0}M$ as Radon measures.  The upper and lower bounds make this measure absolutely continuous with an essentially bounded density.  Local uniformly parabolic Bellman estimates then give \eqref{eq:positive-regularity}; consequently that density is the classical $v_{xx}$ for positive time and the inequalities hold pointwise.  This closes the logical circle without identifying the value function with a classical HJB solution before comparison and stability have identified the viscosity solution.

\subsection{Displacement estimates for the derivative and its inverse}
For a fixed admissible action process write $X_\tau^{x,I}=x+Y_\tau^I$.
Bounded coefficients and Burkholder--Davis--Gundy give
\begin{equation}\label{eq:Y-bound}
 \sup_I\E|Y_\tau^I|
 \le\overline\beta T_0+\sqrt{\overline aT_0}
 =:C_{T_0}.
\end{equation}
Let $I^{\varepsilon,x}$ be $\varepsilon$-optimal at $x$, and abbreviate
$Y=Y_\tau^{I^{\varepsilon,x}}$ and
$D_hv=v(\tau,x+h)-v(\tau,x)$.  The same action process is admissible at
$x+h$, and therefore, for $h>0$,
\[
\begin{aligned}
 D_hv&\le \E\!\left[G(x+h+Y)-G(x+Y)\right]+\varepsilon\\
 &=h\E\int_0^1G'(x+Y+\theta h)\,d\theta+\varepsilon.
\end{aligned}
\]
Using $|G'(z)-G'(x)|\le M|z-x|$ and \eqref{eq:Y-bound} gives
\[
 \frac{v(\tau,x+h)-v(\tau,x)}h-G'(x)
 \le M(C_{T_0}+h)+\frac\varepsilon h.
\]
Using an $\varepsilon$-optimal action process at $x+h$ in the reverse
comparison gives the corresponding lower bound.  First
$\varepsilon\downarrow0$ and then $h\downarrow0$ prove
$|v_x(\tau,x)-G'(x)|\le K_{T_0}$.  If
$p=v_x(\tau,x)=G'(I_0(p))$, then
\[
 m|x-I_0(p)|\le |G'(x)-G'(I_0(p))|\le K_{T_0},
\]
which is the second assertion of \cref{lem:displacement}.

\section{Derivation and uniqueness of the inverse equation}\label{app:inverse}
For $T_0\le T$, write $Q_{T_0}:=(0,T_0)\times\R$ and let
$\mathscr D_{T_0}$ be the class of functions satisfying
\[
 \begin{gathered}
 \Xi-I_0\in L^\infty(Q_{T_0}),\\
 \lambda_{T_0}\le\Xi_p\le\Lambda_{T_0}
       \quad\text{a.e. on }Q_{T_0},\\
 \Xi\in L^2_{\rm loc}(0,T_0;H^1_{\rm loc}(\R)),\\
 \Xi_\tau\in L^2_{\rm loc}(0,T_0;H^{-1}_{\rm loc}(\R)),
 \end{gathered}
\]
with at most linear growth, the weak equation \eqref{eq:X-PDE}, and initial
trace $I_0$ in $L^2_{\rm loc}(\R)$.  This is the whole-line weak class used
in \cref{thm:inverse}.

The existence of \(\Xi\) is inherited from the value function and is not a
second PDE existence problem.  Indeed,
\cref{prop:second-derivative,cor:regularity,lem:displacement} imply that
\(v_x(\tau,\cdot)\) is an
increasing bijection of \(\R\), and hence \(\Xi\) is defined by
\eqref{eq:X-def}.  Identity \eqref{eq:Xi-p-basic} and the global
second-derivative bounds give the slope bounds in \eqref{eq:X-slope}; the
displacement estimate gives \(\Xi-I_0\in L^\infty\).

For smooth Bellman approximations, ordinary differentiation of the Legendre
equation gives
\[
 \Xi_\tau^\varepsilon
 =\partial_p\mathcal M^\varepsilon
   (\tau,p,\Xi_p^\varepsilon).
\]
The locally uniform convergence of \(v_x^\varepsilon\) and the common
positive lower bound for \(v_{xx}^\varepsilon\) imply locally uniform
convergence of the inverse maps.  The local energy bounds used in
\cref{app:second-derivative} give weak convergence of their gradients, and strong
monotonicity of the flux identifies the limit by the standard Minty
argument.  Thus the already defined \(\Xi\) satisfies \eqref{eq:X-PDE}
weakly, has the stated initial trace, and belongs to
\(\mathscr D_{T_0}\).  What remains to prove is uniqueness of this inherited
equation in that class.

For uniqueness, let $\Xi^1,\Xi^2\in\mathscr D_{T_0}$ and
$Z=\Xi^1-\Xi^2$.  Define the secant coefficient
\[
 A(\tau,p)=
 \frac{\mathcal M(\tau,p,\Xi_p^1)-\mathcal M(\tau,p,\Xi_p^2)}
      {\Xi_p^1-\Xi_p^2},
\]
when $\Xi_p^1\ne\Xi_p^2$.  On the diagonal choose any Borel value in
\[
 \left[\frac{\underline a}{2\Lambda_{T_0}^2},
       \frac{\overline a}{2\lambda_{T_0}^2}\right].
\]
This choice does not affect $AZ_p$, because $Z_p=0$ on that set.  Since
every branch has
derivative $a_i/(2r^2)$,
\begin{equation}\label{eq:M-monotone}
 \frac{\underline a}{2\Lambda_{T_0}^2}
 \le A\le
 \frac{\overline a}{2\lambda_{T_0}^2}.
\end{equation}
The difference solves $Z_\tau=\partial_p(AZ_p)$.  Choose
$\chi_R=1$ on $[-R,R]$, supported in $[-2R,2R]$, with
$|\chi_R'|\le C/R$.  Testing by $\chi_R^2Z$ and using Young's inequality
gives
\begin{align*}
 \frac12\frac d{d\tau}\int\chi_R^2Z^2\,dp
 &=-\int A Z_p(\chi_R^2Z)_p\,dp\\
 &\le-\frac{\underline a}{4\Lambda_{T_0}^2}
       \int\chi_R^2|Z_p|^2\,dp\\
 &\quad+\frac C{R^2}
 \int_{R<|p|<2R}Z^2\,dp.
\end{align*}
The bounded-displacement condition makes the last term $O(R^{-1})$.
Integrate in time, use the local $L^2$ initial trace, and let $R\to\infty$.
Then $Z=0$ on every compact cylinder.  Thus all approximation subsequences
have the same limit, and $\Xi$ is precisely the inverse of $v_x$.

At a gradient value where two actions tie,
$r\mapsto\mathcal M(\tau,p,r)$ is generally only Lipschitz.  On the common
slope interval, however, every branch derivative lies between
$\underline a/(2\Lambda_{T_0}^2)$ and
$\overline a/(2\lambda_{T_0}^2)$; the same bounds hold for every secant of
the finite maximum.  Thus the maximum-envelope kink does not affect strong
monotonicity, Caccioppoli estimates, or the Minty identification.  Joint
continuity of $\Xi_p$ follows from the identity
\begin{equation}\label{eq:Xi-gradient-identity}
 \Xi_p(\tau,p)=\frac1{v_{xx}(\tau,\Xi(\tau,p))}.
\end{equation}
For positive time, the concave uniformly parabolic estimates give a jointly
continuous representative of $v_{xx}$ on compact cylinders
\cite{Krylov1980,Wang1992}.  Initial-boundary regularity gives the analogous
statement at $\tau=0$ after localization \cite{Wang1992II}.  To verify that
the nonsmooth Bellman approximation satisfies the required bounds uniformly,
apply the fundamental solution of \eqref{eq:zh} to the initial difference
quotient $h^{-1}\int_y^{y+h}G^{\varepsilon\prime\prime}(r)\,dr$.  The
Gaussian tail and mass estimates for bounded first-order fluxes
\cite{Aronson1967,Aronson1968,LadyzhenskayaEtAl1968} yield
\begin{align}
 &\sup_{|x|\le R}
 \left|\frac{z_h^\varepsilon(\tau,x)}h-G''(x)\right|
 \nonumber\\
 &\quad\le C\omega_{G''}(r+h)
 +C\norm{G''}_\infty\sqrt\tau\nonumber\\
 &\qquad+C\norm{G''}_\infty e^{-r^2/(C\tau)}
 +o_\varepsilon(1).
 \label{eq:approximate-identity}
\end{align}
where $\omega_{G''}$ is the local modulus of continuity.  The last term is
uniform for bounded $h$.  Let first $\varepsilon\downarrow0$ and then
$h\downarrow0$ at fixed positive time.  Thus
\begin{equation}\label{eq:vxx-initial-estimate}
 \sup_{|x|\le R}|v_{xx}(\tau,x)-G''(x)|
 \le C\omega_{G''}(r)+C\sqrt\tau
 +C\exp\{-r^2/(C\tau)\}.
\end{equation}
Letting $\tau\downarrow0$ and then $r\downarrow0$ gives
\[
 \lim_{\tau\downarrow0}\sup_{|x|\le R}
 |v_{xx}(\tau,x)-G''(x)|=0,
 \qquad R<\infty.
\]
The displacement estimate
confines $\Xi(\tau,p)$ to a compact state interval when $p$ is bounded.
Thus \eqref{eq:Xi-gradient-identity} is positive and jointly continuous,
which proves \cref{lem:Xp-continuity}.  Continuity of $\Xi_p$ therefore comes
from the positive-time Bellman estimate and the explicit initial-slice
estimate, rather than from weak continuity of $\Xi$.  We do not claim continuity of
$\Xi_{pp}$ across an action junction.

\StrictOrderAppendix

\section*{Funding}
This work was supported in part by the National Natural Science Foundation of
China (NSFC) under Grant Nos.~12271274 and 12571514.

\section*{Conflict of interest}
The author declares that there is no conflict of interest.

\begingroup
\small
\bibliographystyle{plain}
\bibliography{references}

@book{FlemingSoner2006,
  author    = {Wendell H. Fleming and Halil Mete Soner},
  title     = {Controlled {Markov} Processes and Viscosity Solutions},
  edition   = {2},
  publisher = {Springer},
  address   = {New York, NY, USA},
  year      = {2006},
  doi       = {10.1007/0-387-31071-1}
}

@book{Krylov1980,
  author    = {Nikolai V. Krylov},
  title     = {Controlled Diffusion Processes},
  publisher = {Springer},
  address   = {New York, NY, USA},
  year      = {1980},
  doi       = {10.1007/978-1-4612-0348-8}
}

@article{CrandallIshiiLions1992,
  author  = {Michael G. Crandall and Hitoshi Ishii and Pierre-Louis Lions},
  title   = {User's guide to viscosity solutions of second order partial differential equations},
  journal = {Bull. Amer. Math. Soc.},
  volume  = {27},
  number  = {1},
  pages   = {1--67},
  year    = {1992},
  doi     = {10.1090/S0273-0979-1992-00266-5}
}

@article{Aronson1967,
  author  = {Donald G. Aronson},
  title   = {Bounds for the fundamental solution of a parabolic equation},
  journal = {Bull. Amer. Math. Soc.},
  volume  = {73},
  number  = {6},
  pages   = {890--896},
  year    = {1967},
  doi     = {10.1090/S0002-9904-1967-11830-5}
}

@article{Aronson1968,
  author  = {Donald G. Aronson},
  title   = {Non-negative solutions of linear parabolic equations},
  journal = {Ann. Scuola Norm. Sup. Pisa Cl. Sci.},
  volume  = {22},
  number  = {4},
  pages   = {607--694},
  year    = {1968}
}

@article{LvWuZhang2018,
  author  = {Siyu Lv and Zhen Wu and Qing Zhang},
  title   = {Optimal switching under a hybrid diffusion model and applications to stock trading},
  journal = {Automatica},
  volume  = {94},
  pages   = {361--372},
  year    = {2018},
  doi     = {10.1016/j.automatica.2018.04.048}
}

@book{LadyzhenskayaEtAl1968,
  author    = {Olga A. Ladyzhenskaya and Vsevolod A. Solonnikov and Nina N. Ural'tseva},
  title     = {Linear and Quasilinear Equations of Parabolic Type},
  publisher = {American Mathematical Society},
  address   = {Providence, RI, USA},
  year      = {1968},
  doi       = {10.1090/mmono/023}
}

@article{ElKarouiEtAl1987,
  author  = {{El Karoui}, Nicole and Du Huu Nguyen and Monique Jeanblanc-Picqu{\'e}},
  title   = {Compactification methods in the control of degenerate diffusions: Existence of an optimal control},
  journal = {Stochastics},
  volume  = {20},
  number  = {3},
  pages   = {169--219},
  year    = {1987},
  doi     = {10.1080/17442508708833443}
}

@article{HaussmannLepeltier1990,
  author  = {Ulrich G. Haussmann and Jean-Pierre Lepeltier},
  title   = {On the existence of optimal controls},
  journal = {SIAM J. Control Optim.},
  volume  = {28},
  number  = {4},
  pages   = {851--902},
  year    = {1990},
  doi     = {10.1137/0328049}
}

@article{Rath1977,
  author  = {John H. Rath},
  title   = {The optimal policy for a controlled {Brownian} motion process},
  journal = {SIAM J. Appl. Math.},
  volume  = {32},
  number  = {1},
  pages   = {115--125},
  year    = {1977},
  doi     = {10.1137/0132009}
}

@article{ReisingerForsyth2016,
  author  = {Christoph Reisinger and Peter A. Forsyth},
  title   = {Piecewise constant policy approximations to {Hamilton--Jacobi--Bellman} equations},
  journal = {Appl. Numer. Math.},
  volume  = {103},
  pages   = {27--47},
  year    = {2016},
  doi     = {10.1016/j.apnum.2016.01.001}
}

@article{Carrillo1999,
  author  = {Jos{\'e} Carrillo},
  title   = {Entropy solutions for nonlinear degenerate problems},
  journal = {Arch. Ration. Mech. Anal.},
  volume  = {147},
  number  = {4},
  pages   = {269--361},
  year    = {1999},
  doi     = {10.1007/s002050050152}
}

@article{VandeVate2021,
  author  = {{Vande Vate}, John H.},
  title   = {Average Cost {Brownian} Drift Control with Proportional Changeover Costs},
  journal = {Stoch. Syst.},
  volume  = {11},
  number  = {3},
  pages   = {218--263},
  year    = {2021},
  doi     = {10.1287/stsy.2021.0071}
}

@article{OrmeciVandeVate2024,
  author  = {{Ormeci Matoglu}, Melda and {Vande Vate}, John H.},
  title   = {On the Optimality of Stepwise Policies for Managing Capacity, Inventory and Backorders},
  journal = {Stochastics},
  volume  = {96},
  number  = {6},
  pages   = {1700--1734},
  year    = {2024},
  doi     = {10.1080/17442508.2024.2351148}
}

@book{StroockVaradhan1979,
  author    = {Daniel W. Stroock and S. R. Srinivasa Varadhan},
  title     = {Multidimensional Diffusion Processes},
  publisher = {Springer},
  address   = {Berlin, Germany},
  year      = {1979},
  doi       = {10.1007/978-3-662-12505-8}
}

@book{Lieberman1996,
  author    = {Gary M. Lieberman},
  title     = {Second Order Parabolic Differential Equations},
  publisher = {World Scientific},
  address   = {River Edge, NJ, USA},
  year      = {1996},
  doi       = {10.1142/3302}
}

@article{SaldiYukselLinder2016,
  author  = {Naci Saldi and Serdar Y{\"u}ksel and Tam{\'a}s Linder},
  title   = {Near optimality of quantized policies in stochastic control under weak continuity conditions},
  journal = {J. Math. Anal. Appl.},
  volume  = {435},
  number  = {1},
  pages   = {321--337},
  year    = {2016},
  doi     = {10.1016/j.jmaa.2015.10.008}
}

@book{SaldiLinderYuksel2018,
  author    = {Naci Saldi and Tam{\'a}s Linder and Serdar Y{\"u}ksel},
  title     = {Finite Approximations in Discrete-Time Stochastic Control: Quantized Models and Asymptotic Optimality},
  publisher = {Springer},
  address   = {Cham, Switzerland},
  year      = {2018},
  doi       = {10.1007/978-3-319-79033-6}
}

@incollection{Sheng1982,
  author    = {Diane D. Sheng},
  title     = {A class of problems in the optimal control of diffusions with finitely many controls},
  booktitle = {Advances in Filtering and Optimal Stochastic Control},
  editor    = {Wendell H. Fleming and Luis G. Gorostiza},
  series    = {Lecture Notes in Control and Information Sciences},
  volume    = {42},
  pages     = {353--375},
  publisher = {Springer},
  address   = {Berlin, Germany},
  year      = {1982},
  doi       = {10.1007/BFb0004552}
}

@article{Wang1992,
  author  = {Lihe Wang},
  title   = {On the regularity theory of fully nonlinear parabolic equations: {I}},
  journal = {Comm. Pure Appl. Math.},
  volume  = {45},
  number  = {1},
  pages   = {27--76},
  year    = {1992},
  doi     = {10.1002/cpa.3160450103}
}

@article{LuSerfozo1984,
  author  = {F. V. Lu and Richard F. Serfozo},
  title   = {{M/M/1} Queueing Decision Processes with Monotone Hysteretic Optimal Policies},
  journal = {Operations Research},
  volume  = {32},
  number  = {5},
  pages   = {1116--1132},
  year    = {1984},
  doi     = {10.1287/opre.32.5.1116}
}

@article{GeorgeHarrison2001,
  author  = {Jennifer M. George and J. Michael Harrison},
  title   = {Dynamic Control of a Queue with Adjustable Service Rate},
  journal = {Operations Research},
  volume  = {49},
  number  = {5},
  pages   = {720--731},
  year    = {2001},
  doi     = {10.1287/opre.49.5.720.10605}
}

@article{AtaShneorson2006,
  author  = {Baris Ata and Shiri Shneorson},
  title   = {Dynamic Control of an {M/M/1} Service System with Adjustable Arrival and Service Rates},
  journal = {Management Science},
  volume  = {52},
  number  = {11},
  pages   = {1778--1791},
  year    = {2006},
  doi     = {10.1287/mnsc.1060.0587}
}

@article{BonamiLejeune2009,
  author  = {Pierre Bonami and Miguel A. Lejeune},
  title   = {An Exact Solution Approach for Portfolio Optimization Problems Under Stochastic and Integer Constraints},
  journal = {Operations Research},
  volume  = {57},
  number  = {3},
  pages   = {650--670},
  year    = {2009},
  doi     = {10.1287/opre.1080.0599}
}

@article{BertsimasEtAl1999,
  author  = {Dimitris Bertsimas and Christopher Darnell and Robert Soucy},
  title   = {Portfolio Construction Through Mixed-Integer Programming at {Grantham, Mayo, Van Otterloo and Company}},
  journal = {Interfaces},
  volume  = {29},
  number  = {1},
  pages   = {49--66},
  year    = {1999},
  doi     = {10.1287/inte.29.1.49}
}

@article{SubramanianYang2017,
  author  = {Ajay Subramanian and Baozhong Yang},
  title   = {Optimal Dynamic Risk Taking},
  journal = {Mathematics of Operations Research},
  volume  = {42},
  number  = {3},
  pages   = {599--625},
  year    = {2017},
  doi     = {10.1287/moor.2016.0819}
}

@article{WangZenios2022,
  author  = {Zhengli Wang and Stefanos A. Zenios},
  title   = {New Venture Creation: A Drift--Variance Diffusion Control Model},
  journal = {Operations Research},
  volume  = {70},
  number  = {5},
  pages   = {2982--2997},
  year    = {2022},
  doi     = {10.1287/opre.2021.2171}
}

@article{AltLuckhaus1983,
  author  = {Hans Wilhelm Alt and Stephan Luckhaus},
  title   = {Quasilinear Elliptic--Parabolic Differential Equations},
  journal = {Mathematische Zeitschrift},
  volume  = {183},
  number  = {3},
  pages   = {311--342},
  year    = {1983},
  doi     = {10.1007/BF01176474}
}

@article{Otto1996,
  author  = {Felix Otto},
  title   = {{$L^1$}-Contraction and Uniqueness for Quasilinear Elliptic--Parabolic Equations},
  journal = {Journal of Differential Equations},
  volume  = {131},
  number  = {1},
  pages   = {20--38},
  year    = {1996},
  doi     = {10.1006/jdeq.1996.0155}
}

@article{Wang1992II,
  author  = {Lihe Wang},
  title   = {On the Regularity Theory of Fully Nonlinear Parabolic Equations: {II}},
  journal = {Communications on Pure and Applied Mathematics},
  volume  = {45},
  number  = {2},
  pages   = {141--178},
  year    = {1992},
  doi     = {10.1002/cpa.3160450202}
}

@article{DjehicheHamadenePopier2009,
  author  = {Boualem Djehiche and Sa{\"i}d Hamad\`ene and Alexandre Popier},
  title   = {A Finite Horizon Optimal Multiple Switching Problem},
  journal = {SIAM Journal on Control and Optimization},
  volume  = {48},
  number  = {4},
  pages   = {2751--2770},
  year    = {2009},
  doi     = {10.1137/070697641}
}
\endgroup

\end{document}